# GENERALIZATIONS OF CYCLIC POLYTOPES

Tibor Bisztriczky

Dept. of Mathematics and Statistics, University of Calgary, Calgary, Canada

ABSTRACT. Cyclic polytopes have been studied since at least the early 1900's by Caratheodory and others. A generalization is a construction of a class of polytopes such that the polytopes have some of their properties. The best known example is the class of neighbourly polytopes. Cyclic polytopes have explicit facet structures, important properties and applications in different branches of mathematics. In the past few decades, generalizations of their combinatorial properties have yielded new classes of polytopes that also have explicit facet structures and useful applications. We present an overview of these generalizations along with some applications of the resultant polytopes, and some possible approaches to other generalizations.

## 1. INTRODUCTION

We begin with the following observations from [21], [4] and [5], respectively: " Despite the simplicity of the notion of a polytope, our understanding of what properties a polytope may, or may not, have is severely hampered by the difficulty of producing polytopes having certain desired features.", "The greatest progress in understanding the combinatorial properties of convex polytopes has been made for the special class of simplicial polytopes." and "A problem in studying non-simplicial polytopes is the difficulty of generating examples with a broad range of combinatorial types.".

Well-known examples of non-simplicial polytopes are r-fold d-polytopes P where r is a positive integer and P is, say, a pyramid, or a bipyramid or a prism or a prismoid. It is noteworthy that such P are constructed from lower dimensional polytopes. More general constructions have been introduced by, among others, B. Grunbaum, A. Altshuler, I. Shemer, C. Lee, M. Menzel and A. Padrol.

For a d-polytope $P^*$ and a point $v \notin P^*$, both in $\mathbb{R}^d$, we say that v is **beneath** (**beyond**) a facet F of $P^*$ if $v \notin H$, the affine hull of F, and $P^* \cup \{v\}$ is (is not) contained in the same closed half-space determined by H. The inductive construction, in [21] and [3], of a d-polytope P determined by $P^*$ and v is defined in terms of the relation between v and the facets F of $P^*$. A more refined construction of P in [25] ([22]), called Sewing (Generalized Sewing), is defined by the added requirement of a relation between v and a **flag** (under inclusion, a strictly increasing sequence of faces) of P*.

Regarding the explicit facet structure of a d-polytope P, there are complete descriptions of all facets in terms of vertices if d is small or if P has few vertices or if P is constructed from lower dimensional polytopes. We recall that a cyclic polytope is neighbourly; that is, a **neighbourly** P is a generalized cyclic d-polytope with respect to the number of faces of each dimension, or equivalently, if any k≤d/2 vertices of P determine a face of P. If P is neighbourly, then the explicit facet structure of P is known for all P in the following cases: d=4 and P has at most ten vertices (cf. [1] and [2]); d=6 and P has at most ten vertices (cf. [16]). From [23] and [24], we know the bounds for the number of combinatorial types of neighbourly d-polytopes with eleven vertices and d=4, 6 and 7. It is noteworthy that there are thousands of combinatorial types of neighbourly 4-polytopes (6-polytopes) with eleven vertices, and that the bounds are obtained by Sewing and generalizations of Sewing to oriented matroids (Extended Sewing and Gale Sewing).

In the following; we introduce what may be called the "original" cyclic d-polytopes and consider which of their common properties, say, a "d-polytope" Q may possess. Such properties are expressed as conditions that the "vertices" of Q satisfy, and based upon these conditions, if the set of "facets" of Q is determined; that is, the complete vertex-facet relation that is satisfied by a convex d-polytope in $\mathbb{R}^d$, then Q is a combinatorially constructed generalization. The argument that there is such a Q$\subset \mathbb{R}^d$ is determined by constructions such as noted above.

## 2. CYCLIC POLYTOPES.

We start with the moment curve $\Gamma(\mathbb{I})\subset \mathbb{R}^d$ determined by $\Gamma(t)=(t,t^2,...,t^d)$, d≥3 and t$\epsilon$ $\mathbb{I}$, an open interval in $\mathbb{R}$. We note that $\Gamma(\mathbb{I})$ is simple, $C^\infty$ and of **order d** ($|H\cap \Gamma(\mathbb{I})|$≤d for any hyperplane H$\subset \mathbb{R}^d$ ). Thus, if $t_i<t_j<t_k$ in $\mathbb{I}$ then $\Gamma(t_j)\epsilon\Gamma((t_i,t_k))$ and we say that $\Gamma(t_i)$ **separates** $\Gamma(t_j)$ and $\Gamma(t_k)$ on $\Gamma(\mathbb{I})$. We extend the ordering $t_1 < t_2 <...< t_n$ in $\mathbb{I}$ to an ordering $\Gamma(t_1)< \Gamma(t_2)<...< \Gamma(t_n)$ on $\Gamma(\mathbb{I})$.

Let C(n,d) denote the convex hull of n≥d+2 points $\Gamma(t_s)$ with $t_1<t_2<...<t_n$ in $\mathbb{I}$. Then every facet of C(n,d) is a (d-1)-simplex, $\mathcal{V}$={ $\Gamma(t_1)$, $\Gamma(t_2)$,... $\Gamma(t_n)$ } is the vertex set of C(n,d) and a d element $\mathcal{X} \subset \mathcal{V}$ determines a facet of C(n,d) if, and only if, every two points of $\mathcal{V} \setminus \mathcal{X}$ are separated (on $\Gamma(\mathbb{I})$) by an even number of points of $\mathcal{X}$ ; cf. [20], [21] and [17]. The latter property is called **Gale's Evenness Condition (GEC)**, and $\Gamma(t_1)< \Gamma(t_2)<...< \Gamma(t_n)$ is called a **vertex array** of C(n,d). As noted by Grunbaum in [21]; the construction of C(n,d) uses very few properties of $\Gamma(\mathbb{I})$, and thus, C(n,d) may be constructed via other simple $C^\infty$ curves of order d in $\mathbb{R}^d$.

As background, we introduce our notation and assume familiarity with the basic concepts concerning convex polytopes; cf. [17], [21] and [30].

Let P denote a (convex) d-polytope in $\mathbb{R}^d$. For j=-1,0,1,…,d, let $\mathcal{F}_j(P)$ denote the set of j-faces of P with $\mathcal{V}(P) = \mathcal{F}_0(P)$, $\mathcal{E}(P) = \mathcal{F}_1(P)$ and $\mathcal{F}(P) = \mathcal{F}_{d-1}(P)$, the set of facets of P. The **face lattice** $\mathcal{L}(P)$ of P is the collection of all faces of P ordered by inclusion. The **combinatorial properties** of P are the properties of $\mathcal{L}(P)$, and two d-polytopes are **combinatorially equivalent** (≈) if their face lattices are isomorphic. A **combinatorial construction** of P is a construction of a lattice that is the same as $\mathcal{L}(P)$, and a **realization** of P is the geometric construction of a convex polytope Q⊂ $\mathbb{R}^d$ with the property that Q≈P. For a set Y of points in $\mathbb{R}^d$, let conv(Y) (aff(Y)) denote the convex (affine) hull of Y. If Y={$y_1,y_2,…,y_s$} is finite, we set [$y_1,y_2,…,y_s$] = conv(Y) and <$y_1,y_2,…,y_s$> =aff(Y).

Let P⊂ $\mathbb{R}^d$ be a d-polytope with $\mathcal{V}(P) = \{v_1,v_2,…,v_n\}$, n≥d+2. Then P is **cyclic** if P≈C(n,d); that is, P satisfies GEC with respect to a vertex array ( totally ordered $\mathcal{V}(P)$). As a simplification, we order $\mathcal{V}(P)$ by $v_i < v_k$ if, and only if i<k, and say that $v_j$ separates $v_i$ and $v_k$ in $\mathcal{V}(P)$ if $v_i < v_j < v_k$. Let P be cyclic with $v_1 < v_2 <…<v_n$ Then
- { [$v_1,v_{i-1},v_i$], [$v_k,v_{k+1},v_n$ ]|$3 \leq i \leq n$ and $1 \leq k \leq n-2$} if d=3, and
- {[$v_{i1},v_i,v_j,v_{j+1}$], [$v_1,v_k,v_{k+1},v_n$]|$2 \leq i < j \leq n-2$ and $2 \leq k \leq$n-2 } if d=4.

Accordingly, we anticipate different types of generalizations for odd and even dimensional C(n,d), and seek them among the convex hull of points on simple $C^\infty$ curves Ψ($\mathbb{I}$) in $\mathbb{R}^d$ that have some of the same properties as Γ($\mathbb{I}$). We say that Ψ($\mathbb{I}$) is **convex** if |L∩ Ψ($\mathbb{I}$)|≤2 for any line L⊂ $\mathbb{R}^d$ and Ψ($\mathbb{I}$) ⊂bd(conv Ψ($\mathbb{I}$)). Next, Ψ($\mathbb{I}$) is **ordinary** if each point of Ψ($\mathbb{I}$) has an open neighbourhood of order d in Ψ($\mathbb{I}$). We note that Γ($\mathbb{I}$) is convex and ordinary.

### 3. ODD DIMENSIONS

We start with d=3, (x,y,z)ϵ$\mathbb{R}^3$ and consider the simple $C^\infty$ spherical curves $\Psi_m$= $\Psi_m$($\mathbb{I}$) determined by
$$\Psi_m(t) =(\cos(mt)\sin(t), \sin(mt)\sin(t),\cos(t)), t\epsilon\mathbb{I}=(0,\pi) \text{ and } m\epsilon\mathbb{Z}^+.$$
We note that $\Psi_m$ is convex, |H∩ $\Psi_m$($\mathbb{I}$)|<∞ for any plane H⊂ $\mathbb{R}^3$ and if $\Psi_m(t)\epsilon$H, then there is an open neighbourhood U of t in $\mathbb{I}$ such that H∩ $\Psi_m$(U)={ $\Psi_m(t)$} and $\Psi_m$(U) is contained in a closed half-space determined by H or not. In case of the former (latter), we say that H **supports** (**cuts**) $\Psi_m$ at the point $\Psi_m(t)$.

We note that the horizontal planes H: z=c and -1<c<1 yield that |H∩ $\Psi_m$|=1, and hence, (0,0,1) and (0,0,-1) are boundary points of $\Psi_m$. The vertical planes H, through the z-axis $\mathbb{L}$, yield that |H∩ $\Psi_m$|=m, $\Psi_m$ winds m times about $\mathbb{L}$ ,and H cuts $\Psi_m$ at each point of intersection; cf. [10].

Let $Q_m$⊂ $\mathbb{R}^3$ be a 3-polytope with vertex set $\mathcal{V}_m$ = {$q_1,q_2,…,q_n$}, $q_i$= $\Psi_m(t_i)$ and $t_1<t_2<…<t_n$ in $\mathbb{I}$. Let H⊂ $\mathbb{R}^3$ be a plane such that $\mathbb{L}$ ⊂H, H∩ $\Psi_m$([$t_1,t_n$])=H∩ $\mathcal{V}_m$ = { $p_1,p_2,…,p_u$}, $p_j$= $\Psi_m(s_j)$, $t_1<s_1<s_2…<s_u<t_n$ , $z_j=\cos(s_j)$ and u≥ 4. Then
- $z_u<…<z_2<z_1$ and H∩$Q_m$=[ $p_1,p_2,…,p_u$] is a u-gon with the edges [$p_1,p_2$],[$p_{u-1},p_u$] and [$p_j,p_{j+2}$] for $1 \leq j \leq$u-2,

and

- with $s_0=t_1$ and $s_{u+1}=t_n$, the arcs $\Psi_m((s_{j-1},s_j))$ and $\Psi_m((s_j,s_{j+1}))$ are contained in distinct half-spaces of $\mathbb{R}^3 \setminus H$ for j=1,2,…,u.

Accordingly, if $H \cap Q_m \in \mathcal{F}(Q_m)$ then two distinct elements from $\mathcal{V}_m \setminus (H \cap \mathcal{V}_m)$ are separated in $q_1<q_2<…<q_n$ by an even numbers of $p_j$'s, and $Q_m$ satisfies the necessary part of GEC. We use the preceding two properties of $H \cap Q_m \in \mathcal{F}(Q_m)$ in the construction of a generalization of C(n,3), and then of C(n,d),d≥3. We note from [10] that $\Psi_m$ is an ordinary convex curve, and say that a polytope is **Gale** if it satisfies the necessary part of GEC with respect to some vertex array of P.

A 3-polytope $P \subset \mathbb{R}^3$ is said to be **ordinary** if it has a vertex array $x_0<x_1<…<x_n$, say, such that
(O1)* P is Gale with $x_0<x_1<…<x_n$, and
(O2)* for each facet of P, if $y_1<y_2<…<y_u$ is the (induced) vertex array of F then F is a u-gon with the edges $[y_1,y_2],[y_{u-1},y_u]$ and $[y_j,y_{j+2}]$, j=1,2,…u-2.
Such polytopes are realizable and we refer to [10] for examples.

We say that a d-polytope $P \subset \mathbb{R}^d$, d≥ 3, is **ordinary** if it has a vertex array $x_0<x_1<…<x_n$ such that
(O1) P is Gale with $x_0<x_1<…<x_n$, and
(O2) for each facet of P, if $y_1<y_2<…<y_u$ is the vertex array of F then $G_0,G_1,…,G_u$ are the
(d-2)- faces of F with $G_i=[y_{i-d+2},…,y_{i-1},y_{i+1},…,y_{i+d-2}]$ and the convention that $y_j=y_0$ ($y_j=y_u$) if j<0(j>n).
We note that $G_0=[y_0,y_1,…,y_{d-2}]$, $G_u=[y_{u-d+2},…,y_{u-1},y_u]$ and, as C(n,d) is simplicial, it is ordinary.

REMARK 1. From [10] and [12], we cite the following properties of an ordinary d-polytope $P \subset \mathbb{R}^d$ with (respect to the vertex array) $x_0<x_1<…<x_n$, d≥ 3.
1.1 There is an integer k (the **characteristic char(P)** of P) such that d≤ $k$ ≤ n, and $[x_0,x_i] \in \mathcal{E}(P)$ if, and only if, $[x_{n-i},x_n] \in \mathcal{E}(P)$ if, and only if, i=1,2,…,k.
1.2 P has an explicit facet structure; that is, a complete list of vertex-facet incidences.
1.3 If k=n then P is cyclic with $x_0<x_1<…<x_n$.
1.4 If d is even then P cyclic with $x_0<x_1<…<x_n$.
1.5 If d≥ 5 is odd and k=d, then $\mathcal{F}(P)=\{F_0,F_1,…,F_n\}$ with $F_i=[x_{i-d+1},…,x_{i-1},x_{i+1},…,x_{i+d-1}]$ and, $x_j=x_0$ ($x_j=x_n$) if j<0 (j>n).
We refer to [15] for an inductive construction ( starting with a cyclic (2m+1)-polytope with k+1 vertices) of ordinary (2m+1)-polytopes P with the characteristic k in $\mathbb{R}^{2m+1}$, and note that (O2) and 1.5 indicate that there is a second family of polytopes associated with P.

We say that a d-polytope $M \subset \mathbb{R}^d$, d≥ 2, is a **multiplex** if it has a vertex array $y_0<y_1<…<y_n$ with
$$\mathcal{F}(M)=\{[y_{i-d+1},…,y_{i-1},y_{i+1},…,y_{i+d-1}] | i=0,1,…,n\}$$
and, $y_j=y_0$ ($y_j=y_n$) if j<0 (j>n).

REMARK 2. We refer to [11] for an inductive construction (starting with a d-simplex) and the following properties of a d-multiplex $M \subset \mathbb{R}^d$ with $y_0<y_1<…<y_n$, n≥d+1≥4.
2.1 Each facet of M is a (d-1)-multiplex with the induced vertex array.
2.2 Each vertex figure of M is a (d-1)-multiplex with the induced vertex array.
2.3 M is totally self–dual.

2.4 If d≥5 is odd, then M is ordinary with $y_0<y_1<...<y_n$ and char(M) = d.

Since a d-multiplex is a generalized d-simplex, there is a natural generalization of a simplicial polytope (each proper face is a simplex) to a **multiplicial polytope** (each proper face is a multiplex).

We observe that ordinary polytopes may be characterized as Gale and multiplicial, and it is known that cyclic polytopes may be characterized as Gale and simplicial. In view of 1.4, we have non-cyclic generalizations of C(n,d) if d≥3 is odd. In addition, ordinary polytopes are constructible by the Generalized Sewing method of Lee-Menzel.

## 4. EVEN DIMENSIONS

Let d=2m≥4. Again, we seek the convex hull Q of points $x_0,x_1,...,x_n$ (chosen in order of appearance ) on an oriented simple ordinary convex curve $\Psi$ in $\mathbb{R}^d$ such that
- the convex property of $\Psi$ yields that Q is Gale with $x_0<x_1<...<x_n$, and
- the ordinary property of $\Psi$ yields that Q has non-simplicial facets.

The difficulty is that (with one exception in $\mathbb{R}^4$) there is no explicit model of such a $\Psi$. We do know that the convex hull of k≥d+2 points on an order d subarc of $\Psi$ is a cyclic d-polytope, and that $\Psi$ is the union of such subarcs. Thus, we seek ,say, n≥d+2 points $x_1,x_2,...,x_n$ in $\mathbb{R}^d$ such that
- Q=[ $x_1,x_2,...,x_n$ ] is Gale with $x_0<x_1<...<x_n$, and
- there is an integer k such that d+2≤k≤n, $[x_{i+1},x_{i+2},...,x_{i+k}]$ is a cyclic d-polytope for i=0,1,...,n-k, and $[x_{i+1},x_{i+2},...,x_{i+k},x_{i+k+1}]$ is not cyclic for any 0≤i≤n-k-1.

We call k, the **period** of Q, and say that Q with the latter property is **periodically-cyclic** .

For d=2m≥6 and n≥d+2, there is in [13] an inductive construction of a periodically-cyclic Gale d-polytope $Q^n \subset \mathbb{R}^d$ with n vertices. We start with a cyclic d-polytope $Q^k \subset \mathbb{R}^d$ with $x_1<x_2<...<x_k$, k≥d+2, and assume that for n>r>k, $Q^r=[x_1,x_2,...,x_r]$ is constructed in the prescribed manner. Finally, we choose a point $x_n$ so that

(PC1) $x_n \in <x_1,x_{n-k+1},x_{n-k+2},x_{n-1}>$ and for each facet F of $Q^{n-1}$,

(PC2) <F> does not support $Q^n=[Q^{n-1},x_n]$ in the case F∩[ $x_1,x_{n-k+1},x_{n-1}$]=[$x_1,x_{n-1}$],and

(PC3) F is a facet of $Q^n$ in the case F∩[ $x_1,x_{n-k+1},x_{n-1}$]≠[$x_1,x_{n-1}$] and $[x_1,x_{n-k+1},x_{n-k+2},x_{n-1}] \not\subset F$.

It is noteworthy that this construction yields the explicit facet structure of $Q^n$, and that it does not yield a periodically-cyclic Gale d-polytope if d=4 or if d≥ 5 is odd. In addition, it was observed by M. Bayer that the facets of $Q^n$ also belong to a new family of polytopes.

An e-polytope P⊂ $\mathbb{R}^e$,e≥ 3, is a **braxtope** if it has a vertex array $y_0<y_1<...<y_v$ such that v≥e,
$$\mathcal{F}(P)=\{T_0,T_1,...,T_{v-e+1},E_2,E_3,...E_v\}$$
with $T_i=[y_i,y_{i+1},...,y_{i+e-1}]$, $E_j=[y_0,y_{j-e+2},...,y_{j-1},y_{j+1},...,y_{j+e-2}]$ and $y_t=y_0$ ($y_t=y_v$) if t<0 (t>v).

For e≤ 2, an e-braxtope is an e-simplex. Then a d-polytope Q is **braxial** if each proper face of Q is a braxtope via the ordering induced by a fixed vertex array of Q.

REMARK 3. We refer to [8] for the following properties of an e-braxtope P⊂ $\mathbb{R}^e$ with $y_0<y_1<...<y_v$ and $v \geq e+1 \geq 4$.
3.1 If v≤2e-2 then P is a face of a periodically-cyclic Gale 2m-polytope.
3.2 If v≥2e-1 then there is an inductive construction of P.
3.3 P is braxial.
3.4 The vertex figure of P at $y_0$ is an (e-1)-multiplex with the induced vertex array.

REMARK 4. We refer to [7] for the following properties of a Gale and braxial d-polytope Q⊂ $\mathbb{R}^d$ with $x_0<x_1<...<x_n$, n≥d+1≥4.
4.1 If d is odd then Q is cyclic with $x_0<x_1<...<x_n$.
4.2 If d is even then [$x_0,x_1,...,x_{n-1}$] is a Gale and braxial d-polytope with $x_0<x_1<...<x_{n-1}$.
4.3 If d=2m≥6 then there is an s$\epsilon${0,1,...,n-d+1} such that [$x_j,x_n$]$\epsilon \mathcal{E}$(Q) for s≤ j ≤ n-1.
• If s=1 then Q is cyclic $x_0<x_1<...<x_n$.
• If 2≤ s≤ n-d then Q is periodically-cyclic with the period k=n-s+2.
• If s=n-d+1 then Q is a braxtope with $x_0<x_1<...<x_n$.
4.4 If d=2m≥6 then Q is periodically-cyclic and constructed via (PC1), (PC2) and (PC3).

In view of 4.4, it remains to determine if there are Gale and periodically-cyclic d-polytopes Q if d=4 or, if d=2m≥ 6 and Q is not braxial.
For a generalization of C(n,4), we turn to the generalized trigonometric moment curves Σ and the **bi-cyclic 4-polytopes** B(p,q,n) introduced by Z.Smilansky in [27] and [28]. From [16], we have that the convex hull of n points on the moment curve
$$(\cos 2\pi t, \sin 2\pi t, \cos 4\pi t, \sin 4\pi t), t\epsilon \mathbb{I} = [0,1),$$
is a cyclic 4-polytope. For relatively prime integers q>p>1, Smilansky introduced B(p,q,n) as the convex hull of n evenly spaced points on
$$\Sigma(t) = (\cos 2\pi pt, \sin 2\pi pt, \cos 2\pi qt, \sin 2\pi qt), t\epsilon \mathbb{I};$$
that is, B(p,q,n)=[$b_0,b_1,...,b_{n-1}$] with $b_i=\Sigma(i/n)$.

The curve $\Sigma = \Sigma(\mathbb{I})$ is closed, ordinary and the property of Σ at a point is the same as the property at any other point. Thus, it follows that
[$b_0,b_j,b_k,...,b_s$] $\epsilon \mathcal{F}(B(p,q,n))$ if, and only if, [$b_t,b_{j+t},b_{k+t},...,b_{s+t}$]$\epsilon \mathcal{F}(B(p,q,n))$ for any integer t. We note that any F$\epsilon \mathcal{F}(B(p,q,n))$ is a simplex, or an antiprism over a regular p-gon, or an antiprism over a regular q-gon.

REMARK 5. We refer to [11] for the following properties of B(p,q,n)$\epsilon \mathbb{R}^4$ with n>pq.
5.1 B(p,q,n) is Gale with $b_0<b_1<...<b_{n-1}$ if, and only if, q divides n; furthermore, it has an explicit facet structure.
5.2 There is a number 0<$t_{pq}$<1 such that B(p,q,n) is periodically-cyclic with $b_0<b_1<...<b_{n-1}$ and the period k=[$t_{pq}n$].

The value of $t_{pq}$ is known only for some p and q. For example, $t_{23} \sim 0.419569$ and B(2,3,30) is periodically-cyclic with k=12. In the case that n=pq, q≥7 and q-3≥p≥ $\frac{q+1}{2}$, there are examples of B(p,q,n) that are Gale and periodically-cyclic with the period q.

## 5. APPLICATIONS AND APPROACHES

Two of the important properties of the polytopes listed in the previous sections are that they are non-simplicial and that they have explicit facet structures. Let M(n,d) (O(n,k,d)) denote a d-polytope that is a multiplex (ordinary with the characteristic k). The applicability of M(n,d) and O(n,k,d) to understanding the combinatorial properties of convex polytopes is exemplified below (cf. [17],[21] and the cited articles for definitions).

In [29], R. Stanley observed that the flag vector of M(n,d) is the same as the flag vector of a product of face lattices of polygons , and that the face lattice of M(n,d) is locally self-dual and not the direct product of smaller lattices.

In [19],T. Dinh proved that O(n,k,d) is realizable as a rational polytope, and determined its f-vector for n≥k≥d=2m+1≥5.

In [4], M. Bayer , A. Bruening and J. Stewart determined the flag vector of M(n,d) and, the toric vectors of M(n,d) and O(n,k,5). They showed that the flag vector M(n,d) is the same as the flag vector of a (d-2)-fold pyramid over the (n-d+3)-gon, and presented a computer verification that the set of f-vectors of all O(n,k,d), 5 ≤d=2m+1≤37, spans the Eulerian hyperplane.

In [5], M. Bayer determined the flag vectors of multiplicial d-polytopes for d=2m+1≥5, and the toric vector of O(n,k,d).In addition, there is a construction of multiplicial polytopes with a facet with a large number of vertices.

In [6], M. Bayer presented shallow triangulations of M(n,d) and O(n,k,d), and a combinatorial interpretation of the h-vector of O(n,k,d) based on a shelling of O(n,k,d).

In [9], L. Billera and E. Nevo constructed Bier posets (over the face posets of M(n,d),n≥d≥4) that are nonpolytopal nonsimplicial semi-lattices with nonnegative toric vector.

We recall that ordinary polytopes are defined via their facets (multiplexes), and that O(n,k,d) is cyclic if k=n, or if d=2m. Also, simplicial polytopes are multiplicial, and a 2m-polytope is cyclic if, and only if, it is Gale and simplicial. Thus, ordinary d-polytopes are natural generalizations of cyclic d-polytopes for all d≥3, and their characterization as "Gale and multiplicial" d-polytopes conveys the idea that "Gale and braxial" d-polytopes, say Q*(d), may be also natural generalizations of cyclic d-polytopes for all d≥3. It is not clear if that is the case as Q*(2m+1) is cyclic and, odd-dimensional cyclic polytopes are not necessarily periodically-cyclic. We know also that Q*(2m) is realizable only if it is periodically-cyclic and then, only for m≥3. In that sense, are there better constructions of "cyclic like" even-dimensional polytopes? Specifically, are there constructions of classes $\mathcal{C}= \{Q(d)|d≥4\}$ of d-polytopes with an explicit facet structure and any of the following properties:
- Every $Q(2m)\epsilon\mathcal{C}$ is periodically-cyclic, Gale and not braxial.
- Every $Q(2m)\epsilon\mathcal{C}$ is Gale and not periodically-cyclic .
- Every $Q(2m)\epsilon\mathcal{C}$ is periodically-cyclic and not Gale.

We note from Section 4 that there are periodically-cyclic Gale 4-polytopes that are bi-cyclic and not braxial. An examination of $B(p,q,n)$ indicates that there may be combinatorial constructions of "cyclic like" d-polytopes $P(n,d)$ with $n≥d+3$, $d≥4$ and a vertex array $v_1<v_2<...<v_n$ so that with $P(r,d)=[v_1,v_2,...,v_r]$ and $r=k+1,k+2,...,n$ for some $d+3≤k≤n$:
i) $P(k,d)$ satisfies (GEC) with $v_1<v_2<...<v_k$,
ii) $P(k+1,d)$ is not a simplicial d-polytope,
iii) $P(r,d)$ is Gale with $v_1<v_2<...<v_r$, and
iv) $[v_1,v_2,...,v_{r-1}]≈[v_2,v_3,...,v_r]$ via the lattice isomorphism induced by $v_j \to v_{j+1}$ .

$C(n,2m)$ satisfies i),iii) and iv), and there a different ways for a $P(n,d)$ to satisfy ii). For example, if $d=4$ and $P(k+1,4)$ is not simplicial with simplicial facets $[v_1,v_2,v_3,v_4,v_{k+1}]$ and $[v_1,v_{k-2},v_{k-1},v_k,v_{k+1}]$, then we obtain from [15] that there is a combinatorial construction of $P(n,4)$ such that $\mathcal{L}(P(n,4))$ is the same as $\mathcal{L}(B(p,q,n))$ in the case that $q=k≥7$, $q-3≥p≥(q+1)/2$ and $n=pq$ . If $d>4$, are conditions i)-iv) sufficient to yield the explicit facet structure of a realizable $P(n,d)$ that is Gale and not braxial?

We note that if P is cyclic with $x_0<x_1<...<x_n$ then P is also cyclic with $x_n<x_{n-1}<...<x_1<x_0$ . Accordingly, is there a combinatorial construction of d-polytopes Q that are
i) Gale with $x_0<x_1<...<x_n$, and
ii) for any $F \epsilon\mathcal{F}(Q)$ , $conv(\{x_{n-j}| x_j\epsilon F\})\epsilon\mathcal{F}(Q)$?

A family of d-polytopes $X_r(d),d≥r+2$, with a vertex array $x_0<x_1<...<x_n$ and an explicit facet structure is introduced in [8] with the property that $X_1(d)$ is a braxtope . $X_r(d)$ is called an **(r,d)-braxtope**. With the convention that $X_r(d)$ is a d-simplex for $d≤r+1$, it is conjectured that theorems concerning braxtopes have natural analogues for (r,d)-braxtopes. It is open as to what are the properties of d-polytopes, all of whose facets are (r,d-1)-braxtopes.

Finally, we recall that an edge E of a polytope P is said to be **universal** if [E,v] is a face of P for every vertex v of P. From [25], we have that if P is a neighbourly 2m-polytope with $n≥2m+3$ vertices then P has $u≤n$ universal edges; furthermore, $u=n$ if, and only if, P is cyclic. In [26], Shemer showed that $u≠n-1$, and called P **almost-cyclic** if $u=n-2$. It is now natural to ask if there exist non-simplicial d-polytopes Q that are Gale with $x_0<x_1<...<x_n$ and have u universal edges for, say, $u≥n-2$?